\documentclass[a4paper,12pt]{article}

\usepackage[english]{babel}
\usepackage[utf8x]{inputenc}
\usepackage[T1]{fontenc}
\usepackage{amsthm}

\usepackage{verbatim}
\usepackage[a4paper,top=3cm,bottom=2cm,left=3cm,right=3cm,marginparwidth=1.75cm]{geometry}

\usepackage{graphicx}
\usepackage{float}
\usepackage{amssymb, amscd, latexsym, epsfig}
\usepackage[colorinlistoftodos]{todonotes}
\usepackage[colorlinks=true, allcolors=blue]{hyperref}

\newtheorem{theorem}{Theorem}[section]

\newtheorem{definition}[theorem]{Definition}
\newtheorem{remark}[theorem]{Remark}

\title{Lattice Hydrodynamics}
\author{Dennis Sullivan}
\date{in memory of Jean-Christophe Yoccoz}

\begin{document}
\maketitle

\begin{abstract}
 Using the combinatorics of two interpenetrating face centered cubic lattices together with the part of calculus naturally encoded  in combinatorial  topology, we construct from first principles a lattice model of 3D incompressible hydrodynamics on triply periodic three space. Actually the construction applies to every dimension, but has special duality features in dimension three.
\end{abstract}

\section{Introduction}

\begin{figure}[H]
\centering
\includegraphics[width=.75\textwidth, angle = 90]{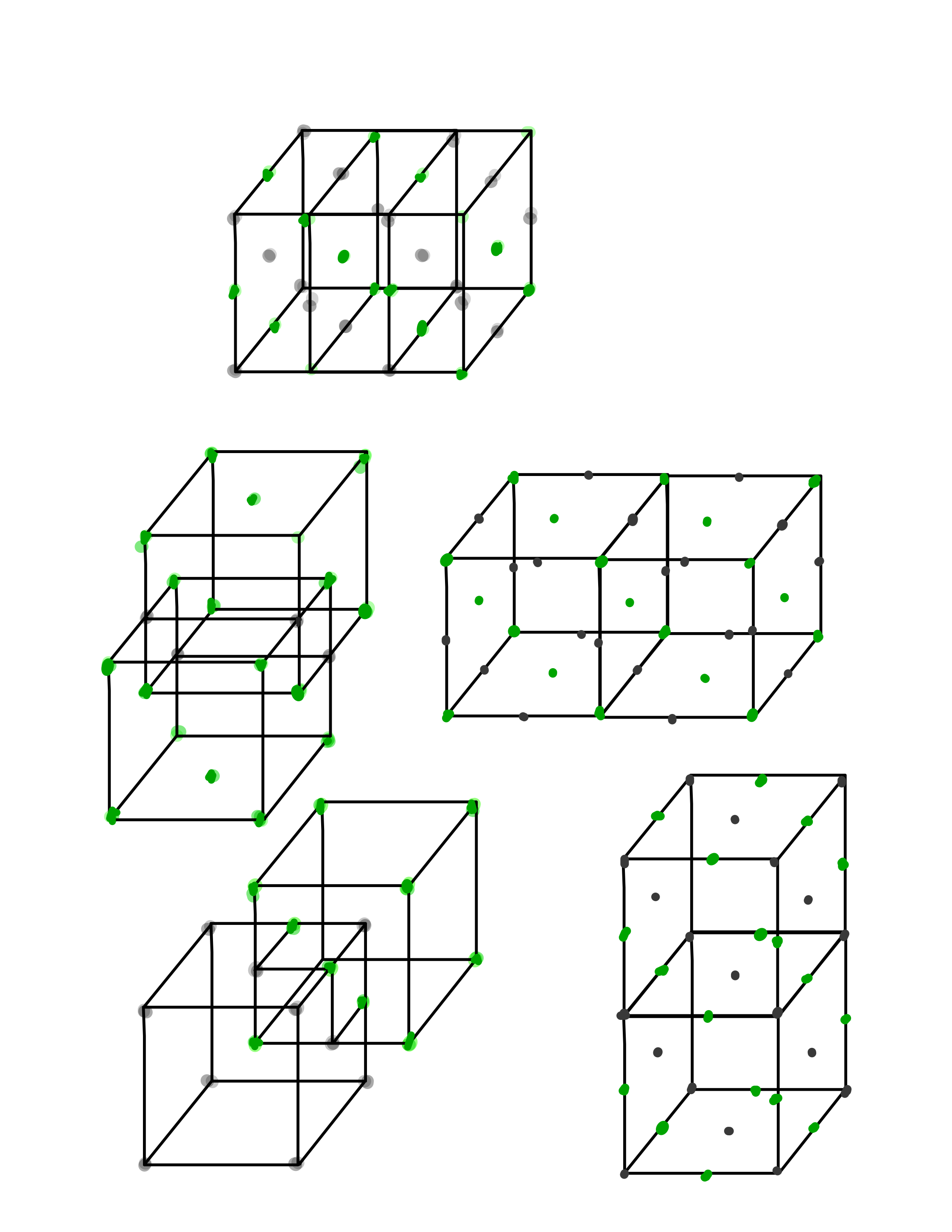}

\caption{\label{fig:fig} How Cubes Intersect}
\end{figure}

We construct a particular lattice model of 3D incompressible fluid motion with viscosity parameter. The construction follows the momentum derivation of  the continuum model switching to combinatorial topology instead of taking the calculus limit.The lattice consists of two interpenetrating face centered cubic lattices which is the crystal  structure of NaCl. The lattice defines sodium  extreme point cubes with their faces, edges  and vertices and chlorine  extreme point cubes with their faces, edges and vertices. In this way the lattice of sites organizes a chain complex $L$ of four vector spaces  built from overlapping uniform cubes, faces, edges and sites giving a multi-layered covering of periodic three space. There are two nilpotent operators on $L$, a duality  involution, each of odd degree, and a combinatorial Laplacian. The result of the momentum derivation is an ODE on one degree of $L$ which is a combinatorial version of  the continuum model.
\begin{equation}\label{eq:eq}
\frac{\partial\{V_L\} }{\partial t} =  \{\ast \delta (\ {V_F} \cdot v_F)\} +  \delta P - \nu \Delta \{ V_L\}, \mbox{ with } \partial\ 
{\{V_L}\}
 =0.
\end{equation}
The combinatorics of the combined lattice $L$ enables a balancing of local and global degrees of freedom required to build the model.

The reference \cite{sullivan} concerns a different approach to models motivated by the  infinite heirarchy of cumulant equations arising from the nonlinearity and its relation to quite modern algebraic topology. There was a difficulty there writing a natural physical model in that context. The model here  is a more recent attempt to  start over from the beginning using  traditional combinatorial topology motivated on the one hand  by observations of fluids and on the other by the classic breakthrough paper of Leray. \cite{leray} 

The goal of work in progress is to use the model both to derive theory and 
to compute meaningfully at a given scale those phenomena that can be naively observed.

\section{The ideas of the construction and definitions}

$L$ denotes the vertices of a regular cubical lattice of edge size $h$ and of even period in three orthogonal directions $(x,y,z)$ which are directed. We imagine a fluid uniformly filling and moving through periodic three space.

\begin{definition} $V_L$: for each site or vertex $q$ of $L$, $V_L(q)$ is a vector at the vertex $q$ which  represents the average velocity of wind or current   taken over the \textbf{cube centered at $q$ with side length $2h$}. Namely the integral times $\frac{1}{8h^3}$. We are assuming the density of particles in the fluid is unity.
\end{definition}

\begin{definition} $V_F, v_F$: For each face $F$of side length $2h$, $V_F$ is  $2h$$V_L(\mbox{center point of $F$})$ and $v_F$ is the component of $V_F$ perpendicular to the oriented $F$ in the direction defined by the right hand rule. 
\end{definition}

\begin{definition} Model proposal:  We are interested for each oriented $F$ in the instantaneous transfer of momentum across the face. This is exactly equal to $\frac{1}{4h^2}$ times the integral over the face $F$ of (the fluid velocity vector times its orthogonal component to $F$). We estimate this average of a product by the product of averages   $V_F \cdot v_F$. This is the closure step in the model which truncates the infinite tower of nonlinear information related to averages of products. When this step proves to be inconvenient the more modern theory of infinity morphisms $[1]$ might be revisited.
\end{definition}

\begin{definition} Since $V_F \cdot v_F$ is a  function on oriented faces of side length $2h$ , we can form  $\delta(V_F \cdot v_F)$, the coboundary of this vector valued function on oriented faces .  This means a vector valued function whose value on an oriented cube (of side length 2h)  is the sum  over its faces of the function on faces,  which are oriented by the outward pointing right hand rule. 
\end{definition}

\begin{definition}   $\ast\delta(V_F \cdot v_F)$  is a lattice vector field, namely a  tangent vector valued function on sites, obtained by placing the value of the  coboundary for the cube at the center of the cube with a sign that depends on the agreement or not of the orientation of the cube with the chosen orientation of space.
\end{definition}

\begin{definition}   $ \{\ast\delta(V_F \cdot v_F)\}:$ The $\{\hspace{.3cm} \}$ of a lattice vector field with $(x,y,z)$ components $(a,b,c)$ at site $q$ is the  one chain obtained by attaching these values to the three edges with center $q$ and length $2h$ in the $(x,y,z)$ directions oriented in their positive sense.This is the bijection between lattice fields and one chains, formalised in the Theorem below but not reiterated there.
\end{definition}

\section{Lattice calculus}

\begin{definition}  Volume preserving: We are modeling fluids that uniformly fill period three space. We say a lattice vector field $V_L(q)$  is \textbf{volume preserving} iff the $1$-chain $\{V_L\}$ from the definition just stated in the previous section has zero boundary, denoted $\partial$. This means if the edges of length $2h$ are re-oriented so the coefficient of $\{V_L\}$ is non negative, then at each vertex the sum of the outgoing coefficients is equal to the sum of the incoming coefficients. This accords with Kirchoff's laws.
\end{definition}

\begin{definition}  Divergence: More generally, the \textbf{divergence}  of a lattice vector field $V_L$ is  $\partial\{V_L\}$ , where  $\{\hspace{.3cm} \}$ is given in the last definition of the previous section. 
\end{definition}

\begin{definition}  Gradient of a lattice scalar field: For a scalar function of vertices $f$ the \textbf{gradient} $f$ is the $1$-cochain whose  value on an oriented edge of length $2h$ is the difference of the values at its two endpoints.  
\end{definition}

\begin{definition} The  Laplacian of a scalar function $f$ of vertices or sites of $L$ is the composition $\Delta f = \partial \delta f$. The value of $\Delta f$ at $q$ is the sum of the values of $f$ at sites $2h$ away from $q$ minus six times the value  of $f$ at $q$. 
\end{definition}

\begin{definition} Curl of a lattice vector field: If $V_L$ is a lattice vector field, 
then  the curl of $V_L$  is the unique lattice vector field that satisfies $\{\mathrm{curl} V_L\}=\ast \delta \{V_L\}$.
\end{definition}
Note: The choice of the edge length $2h$ will be  formalized in the next section. 
\section{Lattice topology, the Laplacian and the Hodge decomposition}\

  For global considerations we need to formalize the choice used above to consider only (and all) positive dimensional cells, i.e. edges, faces and cubes, of side length $2h$. 
    So we consider $L_0$, the vector space generated by the vertices or sites of $L$. Then $L_1$, $L_2$,and $L_3$ are defined respectively to be the vector spaces generated by all the \textbf{oriented} edges, faces and cubes of side length $2h$. This gives twice as many generators as required . This is remedied by imposing the geometric relations (cell, orientation) =-(cell, opposite orientation). Note, as in the figure above, these generators can  overlap. Also at each site  there are exactly three edges of length $2h$ whose midpoint is that site. Thus  dimension $L_1$ = three $\cdot \mbox{ dimension $L_0$}$. This  feature of the choice of side length 2h allows one to confound a lattice vector field with a one chain, which means a linear combination of oriented edges of length $2h$.

\begin{theorem} There are canonical isomorphisms  $\ast:L_0 \leftrightarrow L_3$ and $\ast: L_1  \leftrightarrow L_2$.  If $T$ denotes the tangent space to any point of three space there is a canonical isomorphism $\{\}: L_0 \otimes T  \leftrightarrow L_1$. There are maps $\partial : L_i  \rightarrow L_{i-1}$ and 
$\delta : L_i  \rightarrow L_{i+1}$ satisfying  $\partial \circ \partial=0$, $\delta \circ \delta=0$ and  $\ast \circ \delta=\partial \circ \ast$ , 
 $\ast \circ \partial= \delta \circ \ast$.  Define  $\Delta$ in positive degrees  to be $\partial \circ \delta + \delta \circ \partial$ which extends the previous definition in degree zero.  Then there is an ``orthogonal'' decomposition ( called the decomposition of Hodge) of each $L_i$ as $L_i = \mathrm{im} \partial \oplus \mathrm{im} \delta \oplus \mathrm{kernel} \Delta$.

\end{theorem}

\begin{remark} We note the kernel of $\Delta$ has rank eight in degrees $0$ and $3$ and rank  twenty four in degrees $1$ and $2$. See Note in the Proof. ``Orthogonal'' means relative to the cellular basis, which is orthonormal. 
\end{remark}

\begin{proof}
The graph made of bonds of length $h$ can be two colored because of the even periodicity in all three directions. For a cell of degree one or three of  side length $2h$, there is a center point of one color and $2$ or $8$ vertices in the boundary of the opposite color. For a two cell these corner vertices have the same color as the center point. In general these extreme point vertices of the cells define the vertices of the  cell decomposition of the boundary of the cell used to compute the operators $\partial$ and $\delta$ as is usual in combinatorial topology and Stokes Theorem. Thus a square of side $2h$ has $4$ edges of length $2h$ in its algebraic boundary and a cube of side $2h$ has six faces of edge length $2h$ in its algebraic boundary, etc.

The duality operator $\ast$ relates cells of complementary dimension that intersect transversally at their  center point.The Hodge decomposition is simple and interesting linear algebra valid for any finite dimensional chain complex with positive definite inner product with rational or real coefficients and where the second operator is defined to be the adjoint of the operator defining the chain complex. The kernel of the Laplacian is isomorphic to the homology  [or cohomology] of the complex and defines the "harmonic representatives". Harmonic representatives are both cycles and cocycles, that is, they belong to the intersection of the kernels of the two operators. This follows   in the traditional and interesting way, using the positivity of the inner product after  expanding out ($\Delta$V,V).

The identities are checked pictorially. The signs in the duality isomorphisms are determined by comparing to a global orientation of space.  Note the ordering of dual cells is not important in this comparison because in our odd dimensional space one cell of a dual pair is even dimensional. Otherwise, in even dimensions the order counts half of the time.

Note for the Remark: Since one cells have length 2h  there eight homology classes of vertices. Thus the Laplacian in degree zero has a rank eight kernel.

\end{proof}

\section{The ``potential term'' and the ``friction term''}

The term $\delta P$ in the lattice ODE is meant to cancel the ``volume distortion'' of the ``non linear term'' $\{\ast \delta(V_F \cdot v_F)\}$.  So one wants
$$
-\Delta P = -\partial (\delta P) =  \partial \{\ast \delta(V_F \cdot v_F)\}
$$
In the decomposition of Hodge $\Delta$ preserves the first two factors and is invertible there. Thus we can solve the above and keep the volume preserving property moving forward in time. 

For the last term of the ODE promised above, $-\nu \Delta \{ V_L\}$, one assumes the fluid has a linear response to strain which is isotropic.  This leads in the volume preserving case to  a term proportional to the Laplacian of velocity as explained for example in Landau-Lifschitz "Hydrodynamics". 

Combining all of this we get the ODE equation in words, reading first the LHS and then   the RHS from right to left: ``The rate of change of momentum of a  fluid of uniform unit density inside a cube of side length $2h$ is made up of three parts: 
\begin{enumerate}
\item[i]  the change of momentum due to internal friction, $  \nu \Delta \{ V_L\}. $
\item[ii] the change of momentum $\delta m$ inside the cube created by a potential force of the fluid acting on itself.  The potential $P$ satisfies 
$\{\delta m \}=\delta P$ where $P= \Delta^{-1}(\partial \{\ast \delta(V_F.v_F)\})$.
\item[iii] the change of momentum inside the cube due to a net transfer of momentum across the surface of the cube,  $\{\ast \delta (\ {V_F} \cdot v_F)\}$. 

Thus,

\end{enumerate}

\begin{equation}
 \frac{\partial\{V_L\} }{\partial t} =  \{\ast \delta (\ {V_F} \cdot v_F)\} +  \delta P - \nu \Delta \{ V_L\}, \mbox{ with } \partial\ 
{\{V_L}\}
 =0.
\end{equation}

\end{document}